\documentclass[10pt]{article}
\usepackage[margin=1in]{geometry}
\usepackage{amsmath,amsthm,amssymb,graphicx,verbatim,adjustbox,xcolor}
\usepackage{graphicx,hyperref,tcolorbox,longfbox,framed}
\usepackage{caption}
\newcommand{\bb}[1]{\mathbb{#1}}
\newcommand{\GL}[2]{\text{GL}_{#1}(#2)}

\newtheorem{theorem}{Theorem}[section]

\newtheorem{property}[theorem]{Property}

\newtheorem{algorithm}[theorem]{Algorithm}
\newtheorem{problem}[theorem]{Problem}

\newtheorem{definition}[theorem]{Definition}
\newtheorem{example}[theorem]{Example}
\newtheorem{remark}[theorem]{Remark}
\newtheorem{test}[theorem]{Test}

\title{Conjugacy of Integral Matrices over Algebraic Extensions}

\author{Rebecca Afandi}
\date{\vspace{-5ex}}

\begin{document}
\maketitle

\begin{center}
\textbf{Abstract}
\end{center}

We consider conjugacy of integral matrices by elements in $\GL{n}{R}$ for certain rings $R$ with subring $\bb{Z}$. We note that a Hasse principal does not hold in the context of matrix conjugacy because matrices which are $\GL{n}{\bb{Z}_p}$-conjugate for all $p$ are not necessarily $\GL{n}{\bb{Z}}$-conjugate. By a theorem of Guralnick, we know that integral $n \times n$ matrices are $\GL{n}{\bb{Z}_p}$-conjugate for all primes $p$ if and only if they are conjugate by an element in $\GL{n}{E}$ for some algebraic integral extension $E$ of $\bb{Z}$. We study the problem of finding this extension $E$. Since a result by Latimer and MacDuffee for describing $\bb{Z}$-conjugacy can be generalized to the context of $R$-conjugacy for $R$ any integral domain, we can adapt an existing algorithm for $\bb{Z}$-conjugacy to a new context. We also offer a method for finding $E$ which makes use of the principal ideal theorems of class field theory. We illustrate our method in several examples.

\section{Introduction}
For a field $F$, the theory of conjugacy of matrices in $F^{n \times n}$ by elements in the general linear group, $\GL{n}{F}$, is classical and well-developed. We enter into an area of ongoing research if we work over a ring rather than a field.
For a ring $R$, we say that matrices $A,B \in R^{n \times n}$ are \textbf{$R$-conjugate} or \textbf{conjugate over $R$} to mean that they are conjugate by an element in the general linear group, $\GL{n}{R}$, i.e., there exists $C \in \GL{n}{R}$ with $C^{-1}AC=B$. We denote the equivalence relation given by $R$-conjugacy by $\sim_R$.
 The restriction that the determinant of the conjugating matrix need be a unit in $R$ results in a more nuanced situation in which many open problems emerge. 

In this paper, we consider $R$-conjugacy of integral matrices for $R$ a number of rings, including the ring $\bb{Z}$, algebraic extensions of $\bb{Z}$ (we will denote these by $E$), the ring $\bb{Z}_p$ of $p$-adic integers for a prime $p$, and $\bb{Z}_{(p)}$, the localization of $\bb{Z}$ at a prime $p$. Guralnick \cite{Gur} showed that $\bb{Z}_p$ conjugacy of integral matrices implies conjugacy of those matrices over $\bb{Z}_{(p)}$ (we provide some details of Guralnick's argument in Remark \ref{remark}). Since $\bb{Z}_p$-conjugacy and $\bb{Z}_{(p)}$-conjugacy are equivalent, we will more often refer to $\bb{Z}_p$-conjugacy for ease of notation. Unless specified otherwise, we use $R$ to denote any of the aforementioned rings. 

Much progress has been made in understanding $\bb{Z}$-conjugacy. Grunewald \cite[Thm. A]{Gru} proved that there exists an algorithm for determining whether integral matrices are $\bb{Z}$-conjugate. Grunewald's idea for an algorithm relies on isomorphism-testing of certain modules, however, he did not offer practical means of 
constructing all the necessary submodules. Eick, Hofmann, and O'Brien \cite{EHO} successfully applied Grunewald's ideas 
and implemented an algorithm for $\bb{Z}$-conjugacy.

If we restrict ourselves to integral matrices with a square-free characteristic polynomial, there is an interesting theoretical result which relates the theory of $\bb{Z}$-conjugacy to equivalence of certain number theoretic objects. We will mainly be concerned with this case, so we set the following notation.

 For square-free $f \in \bb{Z}[x]$ of degree $n$, we let $\mathcal{M}_f$ denote the set defined by
\begin{equation*}
\mathcal{M}_f=\{A \in \bb{Z}^{n \times n}: \text{det}(xI-A)=f\}.
\end{equation*}

Note that every matrix in $\mathcal{M}_f$ has the same rational canonical form, the companion matrix of $f$. Thus, the set $\mathcal{M}_f$ gives a single $\bb{Q}$-conjugacy class which we may partition further into $\bb{Z}$-conjugacy classes.

Supposing that $\displaystyle f=\prod_{i=1}^r f_i$ is square-free and each $f_i$ is an irreducible factor with root $\alpha_i$, one of the 
landmark results of the theory of $\bb{Z}$-conjugacy is the LM correspondence \cite{LM}, which states that $\bb{Z}$-conjugacy classes of matrices within $\mathcal{M}_f$ correspond to certain fractional ideal classes in $\displaystyle K:=\prod_{i=1}^r \bb{Q}(\alpha_i)$. Taussky \cite{Tau} and Marseglia \cite{Mar} made the theoretical LM correspondence more concrete by providing bijections in the case that $f$ is irreducible and square-free, respectively. Marseglia used this bijection to implement an algorithm in Magma \cite{Magma} which gives the $\bb{Z}$-conjugacy classes of matrices within $\mathcal{M}_f$ for $f$ square-free.

A class of integral matrices which is of particular relevance in this paper are those which are $\bb{Z}_p$-conjugate for all primes $p$; we say that such matrices are \textbf{locally conjugate}. A special case of a theorem of Guralnick \cite[Thm. 7]{Gur} states that integral matrices are locally conjugate if and only if they are conjugate over an algebraic integral extension of $\bb{Z}$. Note that this means that for integral matrices to be $R$-conjugate for any of the rings of concern to us, we must first have that they are $\bb{Q}$-conjugate. We may therefore express $R$-conjugacy classes as a partition of $\mathcal{M}_f$.

Guralnick's theorem motivates the main problem of this paper.

\begin{problem}\textbf{The conjugacy extension problem}\label{CEP}\\
Let $f \in \bb{Z}[x]$ be square-free. For locally conjugate integral matrices $A$ and $B$ in $\mathcal{M}_f$, find an algebraic extension $E$ of $\bb{Z}$ such that $A$ and $B$ are $E$-conjugate.
\end{problem}

Before discussing our approach to attacking Problem \ref{CEP}, we give a preliminary result regarding $\bb{Z}_p$ conjugacy. 

\begin{theorem}\label{disc}
Let $f \in \bb{Z}[x]$ which is square-free of degree $n$, and take $A, B \in \mathcal{M}_f$. For any prime $p$ which does not divide the discriminant, $\text{disc}(f)$, the matrices $A$ and $B$ are in the same $\GL{n}{\bb{Z}_p}$-conjugacy class.
\end{theorem}

The previous theorem tells us that to determine whether matrices are locally conjugate, one need only test for $\bb{Z}_p$-conjugacy for finitely many primes $p$.

In Section 3.1, we note that the LM correspondence can be generalized to the context of conjugacy over an integral domain. Since the rings $R$ we consider in this paper are integral domains with $\bb{Z}$ as a subring, there is a nice relationship between the original and generalized bijections. If the $\bb{Z}$-conjugacy class of a matrix $A$ corresponds to the fractional ideal class of $I$, then the $R$-conjugacy class of $A$ corresponds to the fractional ideal class of $R \otimes I$. Since a fractional ideal is a free $\bb{Z}$-module, the tensor product essentially extends scalars to $R$. If $f$ factors into the product $\displaystyle \prod_{i=1}^r f_i$ of irreducible factors over $R[x]$, we also identify $R \otimes I$ as being contained within the direct $\displaystyle \prod_{i=1}^r \text{Frac}(R)(\alpha_i)$ by considering $\alpha$ as the tuple $(\alpha_1,..,\alpha_r)$ of the roots of the irreducible factors.

Since the LM correspondences for $R$-conjugacy and $\bb{Z}$-conjugacy can be related via the tensor product, we can adapt an algorithm in \cite{Mar} which tests for $\bb{Z}$-conjugacy to obtain Algorithm \ref{Ralgorithm}, which tests for $R$-conjugacy of integral matrices. We implemented this algorithm in Magma in the case that $R=\mathcal{O}_F$ for a number field $F$ and for matrices with characteristic polynomial which is irreducible in $R[x]$.

The algorithm determines whether a given pair of matrices are $R$-conjugate by determining whether a particular ideal is principal after extension of scalars to $R$. Because of this, it is natural to make use of the principal ideal theorems of class field theory. 
For instance, for a number field $K$, all fractional $\mathcal{O}_K$-ideals are principal in the ring of integers of the Hilbert class field of $K$. 

This leads us to a natural question: ``Does the choice of $E=\mathcal{O}_{K_1}$ for $K_1$ the Hilbert class field of $K:=\bb{Q}(\alpha)\cong \bb{Q}[x]/(f)$ solve the conjugacy extension problem for matrices in $\mathcal{M}_f$?"
We show in Example 3.4 that the answer to this question is `No'. This is because when $f$ has multiple irreducible factors in $R[x]$, the elements of $R \otimes I$ are tuples with number of components equalling the number of irreducible factors of $f$. 
For a fractional ideal $I$, it is not straightforward to test whether $R \otimes I$ is principal unless we are in the case that $f$ is irreducible in $R[x]$. 

We offer a method for trying to solve the conjugacy extension problem which entails finding whether proper subfields of class fields satisfy certain criteria  (see Test \ref{Km} in Section 3.3.1 for details). If such a subfield can be found in a particular example, then its ring of integers gives a solution to the conjugacy extension problem in that instance. We successfully carry out Test \ref{Km} to find an extension $E$ such that matrices in all locally conjugate matrices in $\mathcal{M}_f$ are $E$-conjugate in several examples with irreducible characteristic polynomial.

\subsection*{Acknowledgements}
The author's work has been supported in part by NSF Grant DMS-1720146, which is gratefully acknowledged.

\section{Preliminaries}
Before discussing our contributions to addressing 
the conjugacy extension problem, we discuss what is known about $\bb{Z}$-conjugacy. We also provide some results on $\bb{Z}_p$ conjugacy, including our proof of Theorem \ref{disc}.

\begin{comment}

In this paper, we consider $R$-conjugacy of integral matrices where $R$ denotes any of the following rings:
$\bb{Z}, \bb{Z}_p$ or $\bb{Z}_p$ for a prime $p$, any algebraic extension of $\bb{Z}$. In order for integral matrices to be $R$-conjugate for any of the aforementioned rings $R$, the matrices must be $\bb{Q}$-conjugate. (We will see that $\bb{Z}_p$-conjugacy and $\bb{Z}_p$-conjugacy of integral matrices are equivalent by Remark \ref{remark}, so that $\bb{Z}_p$-conjugate matrices are $\bb{Q}$-conjugate. A theorem of Guralnick \cite{Gur} says that matrices which are conjugate over an algebraic extension of $\bb{Z}$ are $\bb{Z}_{(p)}$-conjugate for every prime $p$ so that they are $\bb{Q}$-conjugate.)

Therefore, it is useful to consider the rational canonical form of a matrix. 

When $f$ is a square-free polynomial in $\bb{Z}[x]$, any matrix with characteristic polynomial $f$ has rational canonical form $\mathcal{C}_f$, the companion matrix of $f$. We focus in this paper on matrices with a given square-free characteristic polynomial; this restrictions allows us to describe matrix conjugacy via a correspondence with certain fractional ideal classes. This motivates the following notation.

We will consider the $R$-conjugacy classes as partitions of the single $\bb{Q}$-conjugacy class given by $\mathcal{M}_f$.
\end{comment}

\subsection{Conjugacy over $\bb{Z}$}
A fundamental result of the study of $\bb{Z}$-conjugacy is a correspondence in \cite{LM}, which we will refer to as the LM correspondence. This correspondence states that the $\bb{Z}$-conjugacy classes of matrices in $\mathcal{M}_f$ with $f$ square-free 
are in bijection with classes of certain fractional ideal classes which are associated to $f$ \cite{LM}.

The LM correspondence was made more explicit in \cite{Tau} in associating the $\bb{Z}$-conjugacy class of a matrix $A$ with irreducible characteristic polynomial with the class of a fractional ideal with $\bb{Z}$-basis comprised of the components of an eigenvector of $A$ \cite{Tau}.

More recently, Marseglia generalized the bijection in \cite{Tau} to the case that $f$ has multiple irreducible factors \cite{Mar}. Actually, the bijection defined in \cite[Thm. 8.1]{Mar} relates matrices with a given square-free minimal polynomial to objects which generalize fractional ideals. We only consider the bijection restricted to the case of square-free characteristic polynomial so that the matrices correspond to fractional ideals, which can be computed using the algorithm in \cite{Mar}. We now state the LM correspondence in terms of Marseglia's formulation.

\begin{theorem}\label{test}(LM correspondence)\cite{LM}, \cite{Mar}\label{LMcor}
Let $\displaystyle f=\prod_{i=1}^r f_i$ be a square-free polynomial in $\bb{Z}[x]$ and let $\alpha=(\alpha_1,...,\alpha_k)$  where $\alpha_i$ denotes a root of the irreducible factor $f_i$.
The $\bb{Z}$-conjugacy classes within $\mathcal{M}_f$ are in bijection with the fractional $\bb{Z}[\alpha]$-ideal classes.\end{theorem}

Following \cite{Mar}, we now discuss what is meant by a fractional $\bb{Z}[\alpha]$-ideal when $\alpha$ denotes a $k$-tuple. If $\alpha=(\alpha_1,..,\alpha_k)$ as in Theorem \ref{LMcor}, then we use $\bb{Z}[\alpha]$ to denote the ring
\begin{equation*}
\bb{Z}[(\alpha_1,..,\alpha_k)]=\{(p(\alpha_1),..,p(\alpha_k)): p \in \bb{Z}[x]\}.
\end{equation*}
The definition of a fractional $\bb{Z}[\alpha]$-ideal in this context is a generalization of the usual notion; a fractional $\bb{Z}[\alpha]$-ideal is defined to be a $\bb{Z}[\alpha]$-submodule of $\displaystyle \prod_{i=1}^r \bb{Q}(\alpha_i)$ which is also a free $\bb{Z}$-module of degree $n=\text{deg}(f)$. Equivalence of $\bb{Z}[\alpha]$-ideals is given by $\bb{Z}[\alpha]$-module isomorphism, which is also the same as scaling by a non-zero-divisor of $\displaystyle \prod_{i=1}^r \bb{Q}(\alpha_i)$.

We denote the set of fractional $\bb{Z}[\alpha]$-ideal classes by $\mathcal{I}_{\bb{Z}[\alpha]}$. Consider the map $\psi$ from the proof of \cite[Thm. 8.1]{Mar}
\begin{align*}
\psi: \mathcal{M}_f &\rightarrow \mathcal{I}_{\bb{Z}[\alpha]}/_{\cong \bb{Z}[\alpha]}\\
A &\mapsto [I]\\
\end{align*}
where $I$ is defined in the following way. For $i=1,..,k$, let $(v_{i1},...,v_{in})^t$ be an eigenvector of $A$ with eigenvalue $\alpha_i$ and 
\begin{equation*}
I=(v_{11},..,v_{k1})\bb{Z} \oplus ...\oplus (v_{1n},...,v_{kn})\bb{Z}.
\end{equation*}
Note that $A$ is the multiplication-by-$\alpha$ map on $I$ with respect to the $\bb{Z}$-basis 
\begin{equation*}
\{(v_{11},..,v_{k1}),...,(v_{1n},...,v_{kn})\}.
\end{equation*}

From $\psi$, we obtain the induced map on $\bb{Z}$-conjugacy classes:
\begin{align*}
\widetilde{\psi}: \mathcal{M}_f/_{\sim \bb{Z}} &\rightarrow \mathcal{I}_{\bb{Z}[\alpha]}/_{\cong \bb{Z}[\alpha]}\\
[A] &\mapsto [\psi(A)].\\
\end{align*}
For a given square-free polynomial $f$ in $\bb{Z}[x]$, Marseglia gives an algorithm to compute all of the fractional $\bb{Z}[\alpha]$-ideal classes \cite{Mar}. This is a new contribution since it was previously only known how to compute the invertible $\bb{Z}[\alpha]$-ideal classes (previous methods are given in \cite{Coh}, \cite{Klu}). Together with the previously discussed bijection, this algorithm can be used to compute a list of representatives of the $\bb{Z}$-conjugacy classes within $\mathcal{M}_f$.

Eick, Hofmann, and O'Brien \cite{EHO} provide an algorithm for testing whether two matrices are $\bb{Z}$-conjugate, even if their minimal polynomial is not square-free. This algorithm relies on isomorphism-testing of certain modules and is motivated by the work of Grunewald in \cite{Gru}.

While the algorithm in \cite{EHO} deals with a more general class of matrices, Marseglia's algorithm ran faster in several examples (as recorded in Table 3 of \cite{Mar}).  Difficulties in efficiently running the algorithm in \cite{EHO} may arise in computing the class group or in testing for isomorphism among a very large number of modules \cite{EHO}. Even the algorithm in \cite{Mar} appears to take exponential time in $n$, the degree of the characteristic polynomial. This is because Marseglia's algorithm entails testing whether a particular fractional ideal is principal, and this is thought to be no easier than computing the class group \cite{Kir}. The most efficient known algorithm is the probabilistic algorithm by Buchmann \cite{Buch}, which is known to have exponential run time in $n$ \cite{Kir}.

With so much recent progress made in understanding $\bb{Z}$-conjugacy, we next shift our focus to $\bb{Z}_p$-conjugacy, and eventually to conjugacy over integral extensions of $\bb{Z}$. 

\subsection{Conjugacy over $\mathbb{Z}_p$}

For a prime $p$, we consider conjugacy of integral matrices by elements in $\GL{n}{\bb{Z}_p}$. This allows us to make some important observations before delving into our discussion of Problem \ref{CEP}.

\begin{remark}\label{remark}
We follow Guralnick's \cite{Gur} argument to show that $\bb{Z}_p$-conjugacy implies conjugacy over $\bb{Z}_{(p)}$. If matrices $A, B \in \mathbb{Z}^{n \times n}$ are conjugate by an element in $\GL{n}{\bb{Z}_p}$, then $A$ and $B$ are $\GL{n}{\bb{Z}/p^k\bb{Z}}$-conjugate for all $k \geq 1$. By the Artin-Rees Lemma, there is a natural number $k'$ such that the $\bb{Z}$-linear map $T(X):=AX-XB$ satisfies the inclusion $p^{k'}\bb{Z}^{n^2} \cap \text{Im}(T) \subseteq p\text{Im}(T)$. The solution $X$ to $AX=XB$ modulo $p^{k'}$ 
may be lifted to $X'$, a matrix in $\GL{n}{\bb{Z}_{(p)}}$ which conjugates $A$ to $B$ (see \cite{Gur}, Theorem 4).
\end{remark}

We now provide a proof for Theorem \ref{disc}.

\begin{proof}
 In the following, we will denote the image of a matrix $A$ modulo $p$ by $A(p)$.

Let $p$ be a prime with $p \nmid \text{disc}(f)$. Since $f$ modulo $p$ is square-free, there is a single $\bb{F}_p$-conjugacy class, containing $A(p)$ and $B(p)$. Let $C \in \GL{n}{\bb{F}_p}$ be a matrix which conjugates $A$ to $B$. We show that $C$ may be lifted to a conjugating matrix in $\GL{n}{\bb{Z}_p}$.

Consider the $\bb{Z}$-linear map $T(X):=AX-XB$. Note that $T(C) \equiv 0$ modulo $p$.

We wish to show that $p\mathbb{Z}^{n^2} \cap \text{Im}(T) \subseteq p \text{Im}(T)$. This is the inclusion following from the Artin-Rees Lemma with $k'=1$. Working with suitable bases for the domain and codomain of $T$, we may express $T$ by a matrix, $S$, in its Smith normal form. Since $f$ is square-free, the rank $r$ of $S$ is $r=n^2-n$. Since also $f$ is square-free modulo $p$, the matrix $S(p)$ also has rank $r$.

Letting $s_1,..,s_r$ denote the non-zero diagonal entries of $S$, we see that
 \begin{align*} p\mathbb{Z} ^{n^2} \cap \text{Im}(T)& = (p\mathbb{Z}  \times ... \times p\mathbb{Z} ) \cap \ (s_1\mathbb{Z}  \times ... \times s_r\mathbb{Z} )\\&=\text{LCM}(p,s_1)\mathbb{Z} \times... \times \text{LCM}(p,s_r)\mathbb{Z}\\&=ps_1\mathbb{Z} \times... \times ps_r\mathbb{Z} \hspace{.3 in} \text{(as } p \nmid s_i)\\
&= p\text{Im}(T).
\end{align*}

By the previous inclusion, we have $T(C)=pT(D)=T(pD)$ for some $D \in \mathbb{Z}^{n^2}$. From this, we obtain $T(C-pD)=T(C)-T(pD)=0$ and $\text{det}(C-pD) \equiv \det{(C)}$, which is non-zero modulo $p$. Thus, $C':=C-pD$ is an element in $\GL{n}{\bb{Z}_{(p)}}$ which conjugates $A$ to $B$.

\end{proof}

We ask whether a local-global principle holds in the context of matrix conjugacy. In other words, we ask whether locally conjugate matrices are conjugate over $\bb{Z}$.

This question is answered affirmatively in \cite{Nebe} for adjacency matrices of Paley and Peisert graphs. In general, however, locally conjugate matrices need not be $\bb{Z}$-conjugate. We use the LM correspondence to obtain an example which demonstrates this.

\begin{example}
Let $f=x^2+5$ and denote a root of $f$ by $\alpha$. The number field $K=\mathbb{Q}(\alpha)$ has class number 2, and $\mathcal{O}_K=\bb{Z}[\alpha]$. A representative of the non-principal fractional ideal class is $I=2 \bb{Z} \oplus (1+\alpha) \bb{Z}$. The trivial ideal class has representative $\bb{Z}[\alpha]$.

Computing the multiplication-by-$\alpha$ matrices with respect to the $\bb{Z}$-bases of each of the ideal class representatives, we find that the matrices $A=\left( \begin{array}{cc} 0 & 1 \\ -5 & 0 \end{array} \right)$ and $B=\left( \begin{array}{cc} -1 & 2 \\ -3 & 1 \end{array} \right)$ are representatives of the two distinct $\bb{Z}$-conjugacy classes within $\mathcal{M}_f$.\\

Note that $\text{disc}(f)=-2^2\cdot 5$. We can easily verify that $A$ and $B$ are $\bb{Z}_p$-conjugate for $p=2$ and $p=5$, which is enough to determine that $A$ and $B$ are locally conjugate by Theorem 2.1.
\end{example}

While locally conjugate matrices are not necessarily $\bb{Z}$-conjugate, a theorem of Guralnick tells us of what can be concluded about such matrices \cite{Gur}. 

\begin{theorem}( \cite{Gur}, Case of Thm. 7)
Integral matrices are locally conjugate if and only if they are conjugate over an algebraic integral extension of $\bb{Z}$.
\end{theorem}

The proof of Theorem 2.4 is based on an existence theorem in \cite{Dade}, which asserts that for a primitive, homogeneous form, there is an algebraic extension of $\bb{Z}$ over which the form realizes a unit. Dade provides a method for finding such an extension, though there is no bound on the degree of the output, as noted in \cite{Wat}.

This leads us to the main problem we address in this paper, the conjugacy extension problem.

\section{The conjugacy extension problem}

For the remainder of the paper, we discuss our contributions to Problem \ref{CEP}, the conjugacy extension problem.

In Section 3.1, we describe how the LM correspondence can be generalized to the context of conjugacy over an integral domain. This can then be used to adapt Marseglia's algorithm to an algorithm which tests for $R$-conjugacy of integral matrices.

In Section 3.2, we introduce how class fields may be used to generate reasonable candidates for the extension $E$ which solves Problem \ref{CEP}. Class fields offer an alternative method to Dade's method in \cite{Dade} for finding such an extension. We will first introduce a special case of the class field method in Section 3.2.2 and then provide it in full generality in Section 3.3.1. A nice feature of our method is that it is motivated by the theory of $R$-conjugacy of integral matrices in that it makes use of the generalized LM correspondence.

\subsection{Generalizing to an integral domain}
The LM correspondence given by Theorem 2.1 may be generalized to describe conjugacy over any integral domain. Estes and Guralnick \cite{EG} give a theoretical proof showing that, for irreducible $f$, the LM correspondence holds over an integral domain for matrices in $\mathcal{M}_f$. The result in \cite{EG} may be generalized to the case that $f$ is square-free with multiple irreducible factors (the details of this generalization are found in \cite{Afandi}).

More concretely, one can see that LM correspondence holds over an integral domain by adapting the bijection $\psi$ from \cite{Mar}. Since there is an embedding of $\bb{Z}$ in each of the rings $R$ we consider in this paper, we may describe the $R$-conjugacy classes of matrices in $\mathcal{M}_f$ as being in correspondence with $R[\alpha]$-module isomorphism classes in $R \otimes_{\bb{Z}} \mathcal{I}_{\bb{Z}[\alpha]}:=\{R \otimes I: I \in \mathcal{I}_{\bb{Z}[\alpha]}\}$. For a fractional $\bb{Z}[\alpha]$-ideal, $I=\oplus v_i \bb{Z}$, we have $R \otimes I=\oplus v_i R$. Suppose that in expressing $I$ as $I=\oplus v_i \bb{Z}$, each $v_i$ is an integral polynomial in $\alpha$ (we may make this assumption about the $v_i$ since we may clear denominators of the $\bb{Z}$-basis of a fractional ideal to obtain an equivalent fractional ideal). Recall the point that when we express $R \otimes I$ as having $R$-basis element $v_i$, we mean the same polynomial in $\alpha$, but where $\alpha$ denotes the tuple of roots of $f$ in $R[x]$. 

We make the connection between the the original and generalized LM bijections more explicit. To distinguish the type of conjugacy class to which we are referring, we denote the $R$-conjugacy class of $A$ by $[A]_R$. If $ \widetilde{\psi}([A]_{\bb{Z}})=[I]$, then one may define the generalized correspondence by:
\begin{align*}
\widetilde{\psi}_R: \mathcal{M}_f/_{\sim R} &\rightarrow R \otimes \mathcal{I}_{\bb{Z}[\alpha]}/_{\cong R[\alpha]}\\
[A]_R &\mapsto [R \otimes_{\bb{Z}} I].
\end{align*}
Since the LM correspondence generalizes in this straightforward way, we may also easily adapt Marseglia's algorithm in \cite{Mar} to testing whether integral matrices are $R$-conjugate for $R$ an integral domain with subring $\bb{Z}$. 

Before giving this modified algorithm, we need a definition.

\begin{definition}\cite{Mar} For $I$, a $\bb{Z}[\alpha]$-ideal in $K$, the \textbf{multiplicator ring of $I$} is given by the colon ideal of $I$ with itself: $(I:I):=\{x \in \mathcal{O}_K: xI \subseteq I\}$. Letting $\hat{I}=R \otimes I$, the multiplicator ring of $\hat{I}$ is $(\hat{I}:\hat{I}):=\{x \in R\otimes \mathcal{O}_K: x\hat{I} \subseteq \hat{I}\}=R \otimes (I:I)$. Note that $(\hat{I}:\hat{I})$ is the largest subring $\mathcal{O}$ of $R \otimes \mathcal{O}_K$ containing $R[\alpha]$ such that $\hat{I}$ is an $\mathcal{O}$-module.
\end{definition}

\noindent 
\begin{framed}
\begin{algorithm}\label{Ralgorithm}
\underline{\textbf{$R$-Conjugacy of Integral Matrices}}(Adaptation of Algorithm in \cite{Mar})
\end{algorithm}

\textit{Input:} Integral matrices $A, B \in \mathcal{M}_f$ for $f$ square-free and an integral domain $R$ containing $\bb{Z}$. Still let $\alpha$ denote the $k$-tuple of roots of the irreducible factors of $f$ in $R[x]$. Let $\displaystyle K=\prod_{i=1}^r \text{Frac}(R)(\alpha_i)$.\\

\textit{Output:} If $A \sim_R B$, returns `Yes' and the conjugating matrix in $\GL{n}{R}$. Otherwise, returns `No'.\\

\textit{Step 1)} Compute fractional $R[\alpha]$-ideals $\hat{I}$ and $\hat{J}$ corresponding to $A$ and $B$. If $\widetilde{\psi}([A]_{\bb{Z}})=[I]$ and $\widetilde{\psi}([B]_{\bb{Z}})=[J]$, then the $R$-conjugacy classes of $A$ and $B$ correspond to the fractional $R[\alpha]$-ideal classes of $\hat{I}=R \otimes_{\bb{Z}} I$ and $\hat{J}=R \otimes_{\bb{Z}}J$, respectively.\\

\textit{Step 2)} Compute $(\hat{I}:\hat{I})$ and $(\hat{J}: \hat{J})$. If $(\hat{I}:\hat{I}) \neq (\hat{J}:\hat{J})$, then $\hat{I}$ is not in the same class as $\hat{J}$, and so $A$ and $B$ are not $R$-conjugate. Output `No' and terminate. Otherwise, set $\mathcal{O}=(\hat{I}:\hat{I})=(\hat{J}:\hat{J})$ and proceed to Step 3).\\

\textit{Step 3)} Test if $(\hat{I}:\hat{J})\mathcal{O}$ is principal. If not, the fractional ideals $\hat{I}$ and $\hat{J}$ are not equivalent and so $A \nsim_R B$. Output `No' and terminate. Otherwise, there is a $\gamma \in K$ such that $(\hat{I}:\hat{J})=(\gamma)$, meaning that $\hat{I}=\gamma \hat{J}$ and $A \sim_R B$. To compute the conjugating matrix, proceed to Step 4).\\

\textit{Step 4)} Suppose that $\hat{I}=\oplus v_iR$ and $\hat{J}=\oplus w_iR$ (we are expressing each of the fractional ideals with respect to an $R$-basis so that $A(v_1,...,v_n)^t=\alpha (v_1,...,v_n)^t$ and $B(w_1,...,w_n)^t=\alpha (w_1,...,w_n)^t$). Compute the change of basis matrix, call it $C$,  between $\{v_1,..,v_n\}$ and $\{\gamma w_1,...,\gamma w_n\}$. Since $C$ gives the change of basis for two $R$-bases of $\hat{I}$, the determinant of $C$ is a unit in $R$, and $C$ conjugates $A$ to $B$. Return `Yes' and the matrix $C$. Terminate.
\end{framed}

In the context of Problem \ref{CEP}, we begin by assuming that we our input matrices $A$ and $B$ are locally conjugate. If $I$ and $J$ correspond to $[A]_{\bb{Z}}$ and $[B]_{\bb{Z}}$ respectively, the generalized LM correspondence tells us that $\bb{Z}_p \otimes I \cong \bb{Z}_p \otimes J$ as $\bb{Z}_p[\alpha]$-modules for any prime $p$. Using more classical language, we say $I$ and $J$ are in the same \textit{genus} (in the sense of \cite[Ch.6]{Rein}). Is is known that $I$ and $J$ are in the same genus if and only if $(I:I)=(J:J)$ (see Proposition 4.1 and Remark 4.3 of \cite{Mar}). Then in Step 2) of Algorithm \ref{Ralgorithm}, we will have $(\hat{I}:\hat{I})=R \otimes (I:I)=R\otimes(J:J)=(\hat{J}:\hat{J})$. While this equality is ensured for locally conjugate matrices, we still compute $\mathcal{O}$ since it is needed for Step 3).

\subsubsection{Implementation of the algorithm}

We implemented Algorithm \ref{Ralgorithm} in Magma \cite{Magma}
in the case that $R$ is the ring of integers of a number field and the input matrices are in $\mathcal{M}_f$ with $f$ irreducible in $R[x]$. To carry out Step 3. of the algorithm, we use the \textsl{IsPrincipal} function in Magma, which in the case of more than one irreducible factor, can only test for principality of fractional ideals in an \'etale $\bb{Q}$-algebra. If $R \supsetneq \bb{Z}$ and we wish to work with $R[\alpha]$-ideals within the $\text{Frac}(R)$-algebra $\displaystyle \prod_{i=1}^r \text{Frac}(R)(\alpha_i)$, the IsPrincipal function only applies in the case that $r=1$.

In the next section, we will consider whether the principal ideal theorems of class field theory can be applied to solving Problem \ref{CEP}. We will see that care must be taken, as complications arise when $f$ is not irreducible. In most of our examples, we choose an algebraic extension of $\bb{Z}$ so that the characteristic polynomial of our matrices remains irreducible in the extension. This restriction allows us to make use of our implementation.

\subsection{Hilbert class fields}

As noted in Algorithm \ref{Ralgorithm}, integral matrices are conjugate over an algebraic extension $R$ of $\bb{Z}$ if a particular fractional $\bb{Z}[\alpha]$-ideal is principal after extending scalars to $R$. For this reason, it is natural to ask whether matrices are conjugate over class fields, which come equipped with principal ideal theorems. In this section, we discuss the Hilbert class field of a number field.

We refer to \cite{child} for the following definitions and properties of class field theory.

\begin{definition}
The \textbf{Hilbert class field} of a number field $K$ is the maximal abelian unramified extension of $K$. We will denote the Hilbert class field of $K$ by $K_1$ (we denote the Hilbert class field this way in order to be consistent with the conventional notation for ray class fields, as we will see in Section 3.3).
\end{definition}

It is known that $[K_1:K]=h_K$, the class number of $K$. We consider Hilbert class fields due to the following property.

\begin{property}\label{HCFprincipal}
For $K$ a number field with Hilbert class field $K_1$, every fractional $\mathcal{O}_K$-ideal is principal in $\mathcal{O}_{K_1}$.
\end{property}

It is important to note that Property \ref{HCFprincipal} holds for fractional $\mathcal{O}_K$-ideals in the classical sense, i.e., fractional ideals which are submodules of a number field, but not necessarily for ideals within a product of number fields.
We therefore restrict ourselves to the case of matrices in $\mathcal{M}_f$ with $f$ irreducible.
 
Let $\alpha$ denote a root of the irreducible polynomial $f$. Then conjugacy of matrices in $\mathcal{M}_f$ coincides with equivalence of fractional $\bb{Z}[\alpha]$-ideals in $K:=\bb{Q}(\alpha)$. Let $R=\mathcal{O}_{K_1}$. Since $\alpha \in \mathcal{O}_{K_1}$, the polynomial $f$ factors further in $R$. Then extending scalars of the fractional ideals of concern to $R$ yields fractional $R[\alpha]$-ideals whose elements have multiple components.
Property \ref{HCFprincipal} of Hilbert class fields may not be helpful in this context since principality in each component does not imply principality of the object as a whole. As mentioned previously, the Magma function IsPrincipal does not apply to these fractional ideals, so computations for determining principality must be done manually. 

\subsubsection{Example}

The next example demonstrates the subtlety that occurs when working with fractional ideals defined within a product of number fields. 

\begin{example}\label{notHCF}
Let $f=x^2+5$ and $K=\bb{Q}(\alpha) \cong \bb{Q}[x]/(f)$. It is known that $h_K=2$, and $I=2\bb{Z} \oplus (1+\alpha)\bb{Z}$ is a fractional ideal which is not principal. From the trivial fractional ideal $\bb{Z}[\alpha]$, we obtain via the LM correspondence the transpose of the companion matrix of $f$. Set $A=\left(\begin{array}{cc} 0 & 1\\ -5 & 0 \end{array} \right)$. The matrix $B=\left(\begin{array}{cc} -1& 2 \\ -3 & 1 \end{array} \right)$ corresponds to $I$. Then we may choose $A$ and $B$ as representatives of the two $\bb{Z}$-conjugacy classes within $\mathcal{M}_f$.

In this example, the full ring of integers of $K$ is $\bb{Z}[\alpha]$. This ensures that any fractional $\bb{Z}[\alpha]$-ideal has multiplicator ring $\bb{Z}[\alpha]$. Since the multiplicator ring of $I$ and $\bb{Z}[\alpha]$ coincide, these fractional ideals are in the same genus, meaning that $A$ and $B$ are locally conjugate matrices.

One can calculate that the Hilbert class field of $K$ is given by $K_1=\bb{Q}(\beta) \cong \bb{Q}[x]/(x^4+12x^2+16)$. Letting $R=\mathcal{O}_{K_1}$, we consider $R$-conjugacy of $A$ and $B$. Since $f$ factors as $f=(x-\alpha)(x+\alpha)$ in $R[x]$, the $R[\alpha]$-ideal of concern is $R \otimes I=(2,2)R \oplus (1+\alpha,1-\alpha)R$. This $R[\alpha]$-ideal is not principal, meaning that the matrices $A$ and $B$ are not $R$-conjugate. This shows that the ring of integers of the Hilbert class field of $\bb{Q}[x]/(f)$ does not automatically solve Problem \ref{CEP} for locally conjugate matrices in $\mathcal{M}_f$.
\end{example}

We will now provide the details for showing that the fractional $R[\alpha]$-ideal $R \otimes I$ from the previous example is not principal.
We let $\overline{\alpha}=(\alpha,-\alpha)$. For any polynomial $p$ in $\alpha$, we use the notation $\overline{p}$ to indicate the same polynomial evaluated in $\overline{\alpha}$. For $r \in R$, we have $\overline{r}=(r,r)$.

We must determine whether 
\begin{equation*}
R \otimes I=\overline{2}R \oplus \overline{(1+\alpha)} R=(2,2)R \oplus (1+\alpha,1-\alpha)R
\end{equation*}
has a generator $(\gamma_1, \gamma_2) \in K_1 \times K_1$. If such a generator exists, then there is a change of $R$-basis  between $\{\overline{2}, \overline{1+\alpha}\}$ and $\{(\gamma_1,\gamma_2), (\gamma_1\alpha,-\gamma_2\alpha)\}$.  Finding whether $R \otimes I$ is principal is equivalent to determining whether there are $\overline{r}_i=(r_i,r_i) \in R$ with 
\begin{align*}
\overline{2}\overline{r_1}+\overline{(1+\alpha)}\overline{r_2}&=(\gamma_1,\gamma_2)\\
\overline{2}\overline{r_3}+\overline{(1+\alpha)}\overline{r_4}&=(\gamma_1\alpha,-\gamma_2\alpha)\\
\end{align*}
and $\overline{r_1r_4}-\overline{r_2r_3} \in R^\times$. Isolating each component, we are searching for a solution $r_i$ to

\begin{equation}
\begin{array}{ccc}
 2r_1+(1+\alpha)r_2 = \gamma_1 & & 2r_3+(1+\alpha)r_4 = \gamma_1\alpha\\
 2r_1+(1-\alpha)r_2 = \gamma_2 & & 2r_3+(1-\alpha)r_4 = -\gamma_2\alpha.\\
  \end{array}
 \end{equation}

When we consider the principality of $R \otimes_{\bb{Z}} I$, it will be useful to work with a $\bb{Z}$-basis $\mathcal{B}$ of $R$. Letting $\beta$, a root of $x^4+12x^2+16$, denote a primitive element of $K_1$, we work with the $\bb{Z}$-basis $\mathcal{B}=\left \{1, \frac{1}{2} \beta, \frac{1}{4}\beta^2,\frac{1}{8}\beta^3 \right\}$ of $R$.
 Denote the $i$-th basis element of $\mathcal{B}$ by $\mathcal{B}_i$.

We know that $I_1:=2R\oplus (1+\alpha)R$ and $I_2:=2R\oplus(1-\alpha)R$ are principal by Property \ref{HCFprincipal}. In fact, both $I_1$ and $I_2$ have a shared generator, call it $\gamma$. If $R \otimes I=(\gamma_1,\gamma_2)$, then we must have $\gamma_1=\gamma u_1$ and $\gamma_2=\gamma u_2$ for $u_i \in R^\times$.
Multiplying through the equations in (1) by $u_1^{-1}$, we see that we are seeking a unit $u$ of $R$ 
so that
\begin{equation*}
\begin{array}{ccc}
 2r_1+(1+\alpha)r_2 = \gamma & & 2r_3+(1+\alpha)r_4 = \gamma \alpha\\
 2r_1+(1-\alpha)r_2 = \gamma u & & 2r_3+(1-\alpha)r_4 = -\gamma u \alpha.\\
  \end{array}
  \end{equation*}
has a solution $r_i$ with $r_1r_4-r_2r_3$ a unit. 

We rewrite this system in terms of the basis $\mathcal{B}$. For $x \in K_1$,
let $M_x$ denote the multiplication-by-$x$ matrix on $R$ with respect to $\mathcal{B}$. For any element $y \in K_1$, let $\overline{y}$ denote the coefficient vector of $y$ in terms of $\mathcal{B}$. Also let $\overline{r}_i=(r_{i1},...,r_{i4})$ denote the the vector of indeterminates. 

Then if we seek a solution to 
\begin{align*}
 2r_1+(1+\alpha)r_2&=\gamma\\
 2r_1+(1-\alpha)r_2&=\gamma u,\\
 \end{align*}
 we can translate this to solving
\begin{align*}
\left(\begin{array}{c|c} M_2 & M_{1+\alpha} \end{array}\right)\left(\begin{array}{c|c} \overline{r}_1 & \overline{r}_2  \end{array}\right)^t&=\overline{\gamma}\\
\left(\begin{array}{c|c} M_2 & M_{1-\alpha} \end{array}\right)\left(\begin{array}{c|c} \overline{r}_1 & \overline{r}_2  \end{array}\right)^t&=M_u\overline{\gamma}\\
 \end{align*}
 for $\overline{r}_1, \overline{r}_2 \in \bb{Z}^4$. Such a solution must satisfy the difference of the equations
 \begin{equation}
 \left(\begin{array}{c|c} 0_{4 \times 4} & M_{2\alpha} \end{array}\right)\left(\begin{array}{c|c} \overline{r}_1 & \overline{r}_2  \end{array}\right)^t=M_{\gamma}\overline{(1-u)}.\\
 \end{equation}
 
 We examine solutions to (2) using the Smith normal form
 of $ \left(\begin{array}{c|c} 0_{4 \times 4} & M_{2\alpha} \end{array}\right)$. There are matrices $P \in \GL{4}{\bb{Z}}$ and $Q \in \GL{8}{\bb{Z}}$ with $P \left(\begin{array}{c|c} 0_{4 \times 4} & M_{2\alpha} \end{array} \right) Q=S$ where
 $S=\left( \begin{array}{cccc|} 2 & & &  \\  & 2 & &  \\ & & 10 & \\ & & & 10  \end{array} \hspace{.1 in} \textbf{\LARGE \scshape 0}_{4 \times 4} \right)$. An integral solution to (2) is $\displaystyle Q \left(\frac{PM_\gamma\overline{(1-u)}}{(2,2,10,10)^t}\right)$ where the division indicated here is componentwise. From the condition that $P M_\gamma\overline{(1-u)}$ must be in $\bb{Z}$-span of $(2,2,10,10)^t$, we find that $u$ must be a unit satisfying $(\overline{1-u}) \in \text{Span}_{\bb{Z}}((-1,-1,1,-4)^t)$.
 
 In other words, there must be a unit of $R$ and $k \in \bb{Z}$ satisfying
 \begin{align*}
1-u=k(-\mathcal{B}_1-\mathcal{B}_2+\mathcal{B}_3-4\mathcal{B}_4), \text{ or }\\
u=1+k (\mathcal{B}_1+\mathcal{B}_2-\mathcal{B}_3+4\mathcal{B}_4)
\end{align*}
 
By a norm argument, we see that $u$ must satisfy
\begin{align*}
\pm1 =N_{K_1/\bb{Q}}(u)
&=\prod_{\sigma \in \text{Gal}(K_1/\bb{Q})} 1+k \sigma(\mathcal{B}_1+\mathcal{B}_2-\mathcal{B}_3+4\mathcal{B}_4)\\
&=500k^4+700k^3+270k^2+10k+1.
\end{align*}

The only integer $k$ which gives a solution is $k=0$, corresponding to $u=1$. There are clearly solutions to 
 \begin{equation*}
 \left(\begin{array}{c|c} 0_{4 \times 4} & M_{2\alpha} \end{array}\right)\left(\begin{array}{c|c} \overline{r}_1 & \overline{r}_2  \end{array}\right)^t=M_{\gamma}\overline{(1-u)}\\
 \end{equation*}
 for $u=1$, and the only solution to 
 \begin{align*}
 2r_1+(\alpha+1)r_2 &=\gamma\\
 2r_1+(-\alpha+1)r_2 &=\gamma
\end{align*}
is $(\gamma/2,0)$.

Setting $u=1$, we repeat the same argument for the system 
\begin{align}
 2r_3+(\alpha+1)r_4 &=\gamma \alpha\\
 2r_3+(-\alpha+1)r_4 &=-\gamma \alpha.
\end{align}

Using the Smith normal form of $ \left(\begin{array}{c|c} 0_{4 \times 4} & M_{2\alpha} \end{array}\right)$ once again, we find that $PM_\gamma (\overline{2\alpha})$ is not divisible by $(2,2,10,10)^t$, meaning that there is no solution to (3). 

We have shown that there is no $u \in R^\times$ so that the system
\begin{equation*}
\begin{array}{ccc}
 2r_1+(1+\alpha)r_2 = \gamma & & 2r_3+(1+\alpha)r_4 = \gamma \alpha\\
 2r_1+(1-\alpha)r_2 = \gamma u & & 2r_3+(1-\alpha)r_4 = -\gamma u \alpha\\
  \end{array}
  \end{equation*}
has a solution. Thus, we have shown that $R \otimes I$ is not principal and, consequently, that $A$ and $B$ are not conjugate over the ring of integers of $K_1$.

\subsubsection{Subfields of the Hilbert class field}

As this example shows, testing whether $R \otimes I$ is principal when $f$ is \textit{not} irreducible in $R[x]$ can be involved, so we restrict our focus to algebraic extensions $R$ of $\bb{Z}$ such that $f$ remains irreducible in $R[x]$.  While $\mathcal{O}_{K_1}$ does not satisfy this criterion, we can still consider proper subfields of the Hilbert class field. This leads us to the following test for searching for a solution to Problem \ref{CEP}.\\

\noindent
\begin{framed}
\begin{test}\label{K1}
\underline{\textbf{Subfields of $K_1$}}
\end{test}
Suppose we begin with $f$ irreducible in $\bb{Z}[x]$ and locally conjugate matrices $A, B \in \mathcal{M}_f$. Let $I$ and $J$ correspond via the LM correspondence to $A$ and $B$, respectively. Since the matrices are locally conjugate, we know that $(I:I)=(J:J)$. Let $\mathcal{O}$ denote the shared multiplicator ring of $I$ and $J$. Let $K=\bb{Q}(\alpha) \cong \bb{Q}[x]/(f)$ and let $K_1$ denote the Hilbert class field of $K$.\\

We search through subfields of $K_1$ and search for a subfield $F$ satisfying the following criteria:
\begin{enumerate}
\item $f$ is irreducible in $\mathcal{O}_F[x]$.
\item $\mathcal{O}_F \otimes (I:J)\mathcal{O}$ is principal
\end{enumerate}

If $F$ satisfies the above criteria, then $\mathcal{O}_F$ is a solution to Problem \ref{CEP} for $A$ and $B$.\\
\end{framed}

This method of searching through the subfields of the Hilbert class field was successful in several examples. The following is a table indicating how often it worked for a particular sample of real quadratic fields, which we obtained from the LMFDB \cite{database}. 

In the first column, we list the irreducible quadratic polynomial $f$. For the sample of number fields $K:=\bb{Q}(\alpha) \cong \bb{Q}[x]/(f)$ represented in the following table, we have that the ring of integers equals $\bb{Z}[\alpha]$. Then the only possible multiplicator ring is $\bb{Z}[\alpha]$,  meaning that all matrices in $\mathcal{M}_f$ are locally conjugate. The fact that $\mathcal{O}_K=\bb{Z}[\alpha]$ also implies that $\mathcal{I}_{\bb{Z}[\alpha]}$ is simply the ideal class group. The third column in the table gives the class number of $K$. 

All the ideal class groups represented in the table are cyclic, and the fourth column lists a matrix $A$ which is chosen to correspond to a generator of the class group. Then $A$ is not $\bb{Z}$-conjugate to the companion matrix of $f$. The last column records whether Test \ref{K1} successfully yields a subfield $F$ of the Hilbert class field satisfying the criteria. If so, the polynomial defining $F$ over $\bb{Q}$ is listed or `Yes' is listed if the polynomial is too lengthy to display. If no suitable subfield of $K_1$ was found, we list `No'. Note that we did not usually check for conjugacy for fields in which $f$ factored further due to the difficulty of testing for principality of fractional ideals.

\begin{table}[h!]
  \begin{center}
  \caption{We use Test \ref{K1} for $\mathcal{M}_f$ with $f$ defining a real quadratic field $K$. The discriminant of the fields $K$ in the table range from from 1 to 100 if $h_K=2$ and from 1 to 500 if $h_K=3$ or $h_K=4$.
}

    \begin{tabular}{ |c|c |c| c | c|} 
   $f$ & disc($f$) & $h_K$ & $A$ &   Test \ref{K1} successful?\\
       \hline
 $x^2-10$ & $2^3\cdot 5$ & 2 &  $\left( \begin{array}{cc} -1 & 3 \\ 3 & 1 \end{array} \right)$  & $x^2-2$  \\
      \hline
 $x^2-15$ & $2^2 \cdot 3 \cdot 5$ & 2 & $\left( \begin{array}{cc} -1 & 2 \\ 7 & 1 \end{array} \right)$ &  $x^2 + 2x - 11$ \\
            \hline
$x^2-x-16$ & $5 \cdot 13$ &  2  &  $\left( \begin{array}{cc} 0 & 2 \\ 8 & 1 \end{array} \right)$ &  $x^2-52$\\
            \hline
      $x^2-x-21$ &  $5 \cdot 17$ &  2 & $\left( \begin{array}{cc} 0 & 3 \\ 7 & 1 \end{array} \right)$ &   $x^2-2205$\\
              \hline
 $x^2-x-57$ & 229 & 3 & $\left( \begin{array}{cc} -2 & 3 \\ 17 & 3 \end{array} \right)$ & $x^3 + 957x^2 + 206910x - 3157132$ \\
\hline
$x^2-x-64$ & 257 & 3 &   $\left( \begin{array}{cc} -1 & 2 \\ 31 & 2 \end{array} \right)$ & $x^3+270x^2-1498824$ \\
\hline
 $x^2-79$ & $2^2\cdot 79$ & 3 & $\left( \begin{array}{cc} -2 & 3 \\ 25 & 2 \end{array} \right)$  & $x^3-66x^2+1089x-1058$ \\
\hline
$  x^2-x-80$ & $3\cdot 107$ & 3  & $\left( \begin{array}{cc} -1 & 2 \\ 39 & 2 \end{array} \right)$ & $x^3-33x+9$ \\
\hline
 $x^2-x-117$ & $7\cdot 67$ & 3 & $\left( \begin{array}{cc} -2 & 3 \\ 37 & 3 \end{array} \right)$ & No\\
\hline
 $x^2-x-118$ & $11\cdot 43$ & 3 &$\left( \begin{array}{cc} 0 & 2 \\ 59 & 1 \end{array} \right)$  & $x^3+90x^2-102168$\\
\hline
$x^2-x-36$ & $5 \cdot 29 $ & 4 & $\left( \begin{array}{cc} -1 & 2 \\ 17 & 2 \end{array} \right)$ & $x^4-44x^2+464$ \\
\hline
 $x^2-82$ & $2^3\cdot 41$  & 4& $\left( \begin{array}{cc} -2 & 3 \\ 26 & 2 \end{array} \right)$ & $x^4-28x^2+32$ \\
\hline
$x^2-x-111$ & $5\cdot 89$ & 4 & $\left( \begin{array}{cc} -2 & 3 \\ 35 & 3 \end{array} \right)$  & Yes \\
\hline
    \end{tabular}
  \end{center}
\end{table}

Because the ideal class groups were cyclic for this sample of number fields, and since $A$ corresponds to the class group generator, the polynomial in the last column defines a number field with ring of integers $E$ satisfying that all matrices in $\mathcal{M}_f$ are $E$-conjugate.

Notice that Test \ref{K1} did not work for the given example of matrices with characteristic polynomial $f=x^2-x-117$. We will find an extension over which matrices in $\mathcal{M}_f$ are conjugate in Example \ref{degree}.

We also give a table recording how often Test \ref{K1} worked for a sample of imaginary quadratic fields defined by the polynomials $f$ in the following table. This sample was again pulled from the LMFDB \cite{database}.

For each number field $K$ defined by $f$ in Table 2, the ring of integers of $K$ is $\bb{Z}[\alpha]$. Once again, all matrices in $\mathcal{M}_f$ are locally conjugate.  Other than the number field defined by $f=x^2+21$, the number fields represented in the table have cyclic ideal class group, and $A$ corresponds to a generator of the group.
\begin{table}[h!]
  \begin{center}
  \caption{We use Test \ref{K1} for $\mathcal{M}_f$ with $f$ defining an imaginary quadratic field $K$. The table lists number fields $K$ with discriminant ranging from -100 to -1 and class number ranging from 2 to 4.}

    \begin{tabular}{ |c|c |c| c | c|} 
 $f$ & disc($f$) & $h_K$  & $A$ & Test \ref{K1} successful?\\
      \hline
 $x^2-x+4$ & $-3\cdot 5$ & 2 &  $\left( \begin{array}{cc} -1 & 2 \\ -3 & 2 \end{array} \right)$ & $x^2+2x+4$ \\
      \hline
          $x^2+5$ & $-2^2\cdot 5$ & 2 & $\left( \begin{array}{cc} -1 & 2 \\ -3 & 1 \end{array} \right)$ & No \\
      \hline
 $x^2+6$ & $-2^3\cdot 3$ & 2 &   $\left( \begin{array}{cc} 0 & 2 \\ -3 & 0 \end{array} \right)$  & $x^2-8x+64$\\
      \hline
$x^2-x+9$ & $-5\cdot 7$ & 2 &  $\left( \begin{array}{cc} -2 & 3 \\ -5 & 3 \end{array} \right)$  &  $x^2+7$\\
      \hline
         $x^2+10$ & $-2^3\cdot 5 $ & 2 & $\left( \begin{array}{cc} 0 & 2 \\ -5 & 0 \end{array} \right)$ &  $x^2+2$\\
      \hline
         $x^2-x+13$ & $-3\cdot 17$ & 2 & $\left( \begin{array}{cc} -1 & 3 \\ -5 & 2 \end{array} \right)$  &  $x^2+8x+19$\\
      \hline
$x^2+13$ & $-2^2\cdot 13$ & 2 &  $\left( \begin{array}{cc} -1 & 2 \\ -7 & 1 \end{array} \right)$ & No\\
      \hline
         $x^2+22$ & $-2^3\cdot 11$ & 2 & $\left( \begin{array}{cc} 0 & 2 \\ -11 & 0 \end{array} \right)$ &  $x^2-40x+576$\\
      \hline
         $x^2-x+23$ & $-7\cdot 13$ & 2 & $\left( \begin{array}{cc} -3 & 5 \\ -7 & 4 \end{array} \right)$ &  $x^2+7 $ \\
      \hline
         $x^2-x+6$ & $-23$ & 3 & $\left( \begin{array}{cc} 0 & 2 \\ -3 & 1 \end{array} \right)$  & $x^3+6x^2+9x-23$ \\
      \hline
 $x^2-x+8$ & $-31$ & 3 & $\left( \begin{array}{cc} -1 & 2 \\ -5 & 2 \end{array} \right)$  & No \\
      \hline
       $x^2-x+15$ & $-59$ & 3 & $\left( \begin{array}{cc} -2 & 3 \\ -7 & 3 \end{array} \right)$  &  $x^3-3x^2-124844$ \\
      \hline
       $x^2-x+21$ & $-83$ & 3 & $\left( \begin{array}{cc} -2 & 3 \\ -9 & 3 \end{array} \right)$  &  $x^3 - 3x^2 - 17107628$\\
      \hline
           $x^2-x+14$ & $-5 \cdot 11$ & 4 & $\left( \begin{array}{cc} 0 & 2 \\ -7 & 1 \end{array} \right)$ & No \\
      \hline
           $x^2+14$ & $-2^3 \cdot 7$ & 4 & $\left( \begin{array}{cc} -2 & 3 \\ -6 & 2 \end{array} \right)$ & No \\
      \hline
           $x^2+17$ & $-2^2 \cdot 17$ & 4 & $\left( \begin{array}{cc} -2 & 3 \\ -7 & 2 \end{array} \right)$ & No \\
      \hline
           $x^2+21$ & $-2^2\cdot 3 \cdot 7 $ & 4  & $\left( \begin{array}{cc} -2 & 5 \\ -5 & 2 \end{array} \right)$  & Yes\\
      \hline
 \end{tabular}
  \end{center}
\end{table}

\newpage

 The following example gives more details regarding the example from the table with non-cyclic class group. 

\begin{example}
For $f=x^2+21$ and $K$ defined by $f$, the Hilbert class field $K_{1}$ has defining polynomial $x^8 + 84x^6 + 3038x^4 -12348x^2 + 405769$.  The matrix $A=\left( \begin{array}{cc} -2 & 5 \\ -5 & 2 \end{array} \right)$ corresponds to the non-principal fractional ideal $I:=5 \bb{Z} \oplus (\alpha+2)\bb{Z}$.

The ideal class group of $K$ is not cyclic, so we must also consider another non-principal $\bb{Z}[\alpha]$-ideal $J:=-7\bb{Z} \oplus \alpha \bb{Z}$.  The matrix corresponding to $J$ is $B=\left( \begin{array}{cc} 0 & -7 \\ 3 & 0 \end{array} \right)$. 

Let $F$ denote the subfield of $K_1$ defined by $x^4-32x^3+1616x^2-21760x+396544$ and let $E=\mathcal{O}_F$. We have that $f$ remains irreducible in $E[x]$ and both $E \otimes I$ and $E \otimes J$ are principal. Then $\mathcal{C}_f \sim_E A \sim_E B$.
\end{example}

When Test \ref{K1} was successful for the sample of polynomials $f$ in Table 2, the resulting extension $E$ is such that  there is a single $E$-conjugacy class within $\mathcal{M}_f$.

Note that for $f=x^2+5$, we already showed in Example \ref{notHCF} that $A$ and $\mathcal{C}_f$ are not conjugate over the Hilbert class field, so they cannot be conjugate over any proper subfield of $K_1$ either. There are several other instances in which this method did not work. This is not surprising since there is no obvious reason for why there must be a subfield of the Hilbert class field satisfying the criteria of Test \ref{K1}. However, Hilbert class fields provide a natural place to look for candidates for a solution to Problem \ref{CEP}.

We list some more data for non-quadratic fields in the appendix.

\subsection{Ray class fields}

We will extend the work of the previous section to ray class fields. Ray class fields are generalizations of the Hilbert class field, differing only in that ray class fields can be ramified at a finite number of primes. Since the condition on ray class fields is not as stringent as requiring \textit{all} primes to be unramified, ray class fields can be large in degree and have many subfields in comparison to Hilbert class fields. Ray class fields also have a principal ideal property. This allows us to adapt Test \ref{K1} to ray class fields, which provides more potential solutions to Problem \ref{CEP}.

We will define ray class fields and list their principal ideal property. More details for the description below can be found in \cite{child} and \cite{milne}.

For a number field $K$, a \textbf{$K$-modulus}, denoted by $\mathfrak{m}$, is the formal product $\mathfrak{m}=\mathfrak{m}_0\mathfrak{m}_\infty$ where $\mathfrak{m}_0$ is an ideal of $\mathcal{O}_K$ and $\mathfrak{m}_\infty$ is a formal product of infinite primes, or real embeddings of $K$.

For a given modulus $\mathfrak{m}$, there exists an associated field satisfying certain properties called the ray class field of $K$ with modulus $\mathfrak{m}$, conventionally denoted $K_\mathfrak{m}$. The modulus is said  to be admissible for $K_\mathfrak{m}$. There exists a least admissible modulus of a class field $K_\mathfrak{m}$ is called the \textbf{conductor} of $K_\mathfrak{m}$ and is denoted by $\mathfrak{f}$.

Given a number field $K$, we may describe the \textbf{ray class field}, $K_\mathfrak{m}$, as the maximal, abelian extension of $K$ which is ramified exactly at the primes which divide the conductor $\mathfrak{f}$ of $K_\mathfrak{m}$.

Note that when $\mathfrak{m}=1$, the ray class field must be unramified at all primes, and so we obtain that $K_1$ is the Hilbert class field of $K$.

In the context of ray class fields, we have the following principal ideal property.

\begin{property}\label{RCFprincipal}\cite{Tan}
Let $K$ be a number field and $\mathfrak{m}$ a $K$-modulus. Suppose the ray class field $K_\mathfrak{m}$ has conductor $\mathfrak{f}=\mathfrak{f}_0\mathfrak{f}_\infty$. Then every fractional 
$\mathcal{O}_K$-ideal relatively prime to $\mathfrak{f}_0$ is principal in $\mathcal{O}_{K_\mathfrak{m}}$.
\end{property}

Due to the previous property, it is natural to extend our work from Section 3.2.2 to ray class fields.

\subsubsection{Subfields of the ray class field}

Despite Property \ref{RCFprincipal}, the fact that $f$ will factor further in $K_\mathfrak{m}$ means that  $\mathcal{O}_{K_\mathfrak{m}}$ does not necessarily solve Problem \ref{CEP}. However, we can adapt Test \ref{K1} to ray class fields and again search for a suitable subfield.

\noindent
\begin{framed}
\begin{test}\label{Km}
\underline{\textbf{Subfields of $K_\mathfrak{m}$}}
\end{test}

The setup is identical to that in Test \ref{K1}. 

Pick a $K$-modulus $\mathfrak{m}$ which is relatively prime to $(I:J)\mathcal{O}$. Let $K_\mathfrak{m}$ denote the ray class field of $K$ with modulus $\mathfrak{m}$.

We search through subfields of $K_\mathfrak{m}$ and search for a subfield $F$ satisfying the same criteria as in Test \ref{K1}.

If $F$ satisfies the above criteria, then $\mathcal{O}_F$ is a solution to Problem \ref{CEP} for $A$ and $B$.\\
\end{framed}

Note that we obtain Test \ref{K1} from Test \ref{Km} by setting $\mathfrak{m}=1$.

Referring to Example \ref{notHCF}, we showed that locally conjugate matrices $A$ and $B$ in $\mathcal{M}_f$ need not be $\mathcal{O}_{K_1}-$conjugate for $K=\bb{Q}(\alpha) \cong \bb{Q}[x]/(f)$. We show that a subfield of a particular ray class field solves Problem \ref{CEP} for the matrices from this example.

\begin{example}
The matrices $A=\left( \begin{array}{cc} 0 & 1 \\ 15 & 0 \end{array} \right)$ and $B=\left( \begin{array}{cc} -1 & 2 \\ -3 & 1 \end{array} \right)$ are representatives of the two $\bb{Z}$-conjugacy classes within $\mathcal{M}_f$ for $f=x^2+5$. Let $K=\bb{Q}(\alpha)\cong \bb{Q}[x]/(f)$. 

Corresponding to these matrices are the fractional $\bb{Z}[\alpha]$-ideals $\bb{Z}[\alpha]=1\bb{Z} \oplus \alpha \bb{Z}$ and $I=2\bb{Z} \oplus (1+\alpha)\bb{Z}$. As previously noted, both of these fractional ideals have multiplicator ring $\bb{Z}[\alpha]$.

We wish to find an algebraic extension $R$ of $\bb{Z}$ in which $R \otimes (I:\bb{Z}[\alpha])=R \otimes I$ is principal. We compute that $N_{K/\bb{Q}}(I)=2$, meaning that we must select a modulus relatively prime to $2\mathcal{O}_K$. Letting $\mathfrak{m}=3\mathcal{O}_K$, we find that $K_\mathfrak{m}$ is the number field defined by $x^8+12x^6+158x^4-228x^2+3721$ and has conductor $\mathfrak{m}$. 

Applying Test \ref{Km}, we find a subfield $F$ of $K_\mathfrak{m}$ given by 
\begin{equation*}
F=\bb{Q}(\beta)\cong\bb{Q}[x]/(x^4-12x^3+158x^2+228x+3721)
\end{equation*}
satisfying the method's criteria. Thus, $A$ and $B$ are $R$-conjugate for $R=\mathcal{O}_F$.

We computed the conjugating matrix in this example. Working with respect to the $\bb{Z}$-basis
\begin{equation*}
\mathcal{B}=\left\{1, \frac{1}{8}(\beta-1), \frac{1}{64}(\beta^2-2\beta+1), \frac{1}{1024}(\beta^3-3\beta^2+3\beta-513)\right\}
\end{equation*}
of $R$, and letting $\mathcal{B}_i$ denote the $i$-th basis element, we find that
\begin{equation*}
C=\left( \begin{array}{cc} -\mathcal{B}_2 & -1-\mathcal{B}_4 \\  3+\mathcal{B}_2+3\mathcal{B}_4& -1-2\mathcal{B}-2-\mathcal{B}_4 \end{array} \right)
\end{equation*}
is an element in $\GL{2}{R}$ which conjugates $A$ to $B$. 
\end{example}

Based on the data that we collected, Test \ref{Km} seems to typically result in extensions of smaller degree than Dade's method from \cite{Dade}. In particular, Dade's method always yields an extension of degree strictly greater than $n$, while we found extensions of degree $n$ in several instances, e.g., see Tables 1 and 2 (in these tables, the degree of extensions actually coincides with the class number of $K$). In the following example, we found a noteworthy
difference in the degree of the extension resulting from each method.

\begin{example}\label{degree}
The matrices $A=\left( \begin{array}{cc} -2 & 37 \\ 3 & 3 \end{array} \right)$ and $B=\left( \begin{array}{cc} 0 & 117 \\ 1 & 1 \end{array} \right)$ are locally conjugate matrices which are not $\bb{Z}$-conjugate. These matrices have characteristic polynomial $f=x^2-x-117$, and $\text{disc}(f)=7 \cdot 67$.
The rational matrices $C_1=\left( \begin{array}{cc} 3 & -6 \\ 0 & 9 \end{array} \right)$ and $C_2=\left( \begin{array}{cc} -10 & -538 \\ 14 & 12 \end{array} \right)$ conjugate $A$ to $B$ and have relatively prime determinants.

Following Dade's method, we find that the quadratic form $f(x,y)=\text{det}(xC_1+yC_2)$ realizes a unit over an algebraic extension of degree larger than 462. We did not explicitly work out the polynomial defining the extension, nor did we solve for $x$ and $y$, since this degree is so large.

On the other hand, Test \ref{Km} requires that we find an algebraic extension $R$ so that $R \otimes I$ is principal for $I=3\bb{Z} \oplus (2+\alpha)\bb{Z}$. Letting $K=\bb{Q}(\alpha)\cong \bb{Q}[x]/(f)$, we pick the modulus to be the $\mathcal{O}_K$-ideal $\mathfrak{p}=7\bb{Z} \oplus (3+\alpha)\bb{Z}$. This is a prime ideal of $\mathcal{O}_K$ with $\mathfrak{p}^2=(7)$.(We do not have results which narrow down which moduli relatively prime to the norm of the ideal should be considered, but it is interesting to note that the modulus is a prime dividing $7$, which divides the discriminant of $f$). 

Letting $\mathfrak{m}=\mathfrak{p}$, we find that the ray class field $K_\mathfrak{m}$ has degree 18 over $\bb{Q}$ and has conductor $\mathfrak{p}$. Searching through the proper subfields of $K_\mathfrak{m}$, we find that the subfield $F$ defined by the cubic polynomial 
\begin{align*}
x^3 + &c_2x^2 x + c_1x+c_0 \text{ with }\\
c_2&=22427531465691\\
c_1&=87019205503941567942935016 \text{ and }\\
c_0&=-169863356476213700999189634845323984727
\end{align*}
satisfies the criteria of Test \ref{Km}. Therefore, the matrices $A$ and $B$ are $\GL{2}{\mathcal{O}_F}$-conjugate, and $\mathcal{O}_F$ is an algebraic extension of degree $3$ over $\bb{Z}$.
\end{example}

In looking for a solution to Problem \ref{CEP}, this method opens up for consideration subfields of the ray class field of any modulus that is relatively prime to $(I:J)\mathcal{O}$. However, the degree and quantity of subfields of ray class fields can be very large. This presents computational difficulties, as testing whether an ideal is principal in a number field of large degree is costly. 

Due to the expense of these computations, we only used Test \ref{Km} with $\mathfrak{m} \neq 1$ for solving Problem \ref{CEP} in a limited number of examples in which the degree of the ray class field was not too large. These examples are listed in the appendix.

\section{Conclusion and Open Problems}

We have seen that Hilbert class fields do not always provide a solution to Problem \ref{CEP}. While we offered Test \ref{Km} which incorporates ray class fields, we could only carry out this method in a few cases due to computational difficulties. Thus,  several open problems remain in this area of research, including the following.

\begin{enumerate}
\item Is there always a modulus $\mathfrak{m}$ for which there is a subfield meeting the criteria of Test \ref{Km}? Is there a way to characterize those moduli for which this holds? This would be a helpful contribution since as of now we only require that the modulus be relatively prime to a particular fractional ideal. In Example 3.7, our choice of modulus, which was related to the discriminant of the characteristic polynomial, proved successful. Should this hold in general?
\item Can one make the $R$-conjugacy algorithm more efficient? The algorithm which tests for principality of ideals seems to take exponential time in $n$, and it is not clear whether this can be improved.
\item Can one write an algorithm which tests whether a fractional ideal defined in an $R$-algebra of the form $\displaystyle \prod_{i=1}^r \text{Frac}(R)(\alpha_i), k>1$ is principal? This would allow us to extend the $R$-conjugacy algorithm to matrices which have square-free characteristic polynomial with multiple irreducible factors.
\end{enumerate}

\section{Appendix}

Throughout the appendix, we obtained all number fields of consideration from the LMFDB \cite{database}.

\subsection{Test \ref{K1} for cubic and quartic number fields}

For the six cubic fields with class number 2 and discriminant ranging from -1,000 to 1,000 from the LMFDB, Test \ref{K1} does not yield a solution to Problem \ref{CEP}. We do not list these fields in a table.

The next table lists cubic polynomials $f$ which give number fields $K$ with discriminant ranging from -1,000 to 1,000 and class number 3. 

In the examples for which Method 3.6 was successful, a polynomial is listed in the last column which defines a number field $F$ such that the number of $\GL{3}{\mathcal{O}_F}$-conjugacy class of matrices in $\mathcal{M}_f$ coincides with the number of local conjugacy classes. We also recorded when Test \ref{K1} failed.
\begin{table}[h!]
  \begin{center}
  \caption{Hilbert class field method for $\mathcal{M}_f$ with $f$ cubic.}
    \begin{tabular}{|c|c |c|c|c|}
 $f$ & disc($f$) & $h_K$ & \#Local conjugacy classes & Test \ref{K1} successful?\\
      \hline
     $x^3-x^2+5x+1$ & $-2^2 \cdot 3 \cdot 7^2$ & 3  & 1 & $x^3-3x^2-60x+251$\\
     \hline
    $ x^3-3x-10$ & $-2^3 \cdot 3^4$ & 3 & 2 & No\\
     \hline
     $x^3-x^2-4x+12$ & $-2^2 \cdot 13^2$ & 3 & 2 & $x^3-3x^2-114x+467$\\
     \hline
     $x^3+6x-1$ & $3^4 \cdot 11$ & 3 & 1 & No\\
     \hline
     $x^3-x^2+5x-6$ & $-7^2 \cdot 19$ & 3 & 1 & $x^3+3x^2-8376x-303407$\\
    \hline
    $x^3-x^2+5x-13$ & $-2^2 \cdot 5 \cdot 7^2$ & 3 & 2 & $x^3+3x^2-60x+127$\\
     \hline
 \end{tabular}
  \end{center}
\end{table}

Test \ref{K1} failed for the seven cubic fields with discriminant ranging from -2,000 to 2,000 and class number 4 or 5 listed in the LMFDB. 

We checked the seven quartic fields with class number 2 and discriminants ranging from -2,500 to 2,500, the four fields with class number 3 and discriminants in -5,000 to 5,000, and the four fields with class number 4 and discriminants from -10,000 to 10,000 from the LMFDB. Of these, only the quartic polynomials $f$ for which Test \ref{K1} was successful are given in the next table. 

The last column gives a polynomial which defines a number field $F$ so that the number of $\GL{4}{\mathcal{O}_F}$-conjugacy classes within $\mathcal{M}_f$ is given by the number of local conjugacy classes.
\begin{table}[h!]
  \begin{center}
  \caption{Hilbert class field method for $\mathcal{M}_f$ with $f$ quartic.}
  \small
    \begin{tabular}{|c|c |c|c|c|} 
 $f$ & disc($f$) & $h_K$ & \#Local conjugacy classes & Test \ref{K1} successful?\\
      \hline
     $x^4+4x^2+1$ & $2^8 \cdot 3^2$ & 2 & 1 & $x^2-10x+73$\\
     \hline
     $x^4+9$ & $2^8 \cdot 3^2$ & 2 & 3 & $x^2-10x+73$\\
     \hline
     $x^4-x^3+4x^2+x+1$ & $2^3 \cdot 23^2$ & 3 & 1 & $x^3+30x^2-45x-12501$\\
     \hline
     $x^4-2x^3+4x^2+2x+1$ & $2^6 \cdot 5^3 $ & 4 & 1 & $x^4+6x^3+111x^2+526x+761$\\
     \hline
     $x^4-x^3+x^2-6x+6$ & $2^3 \cdot 3^2 \cdot 5^3$ & 4 & 1 & $x^4+32x^2+544x^2+4608x+15616$\\
     \hline
     $x^4+5x^2+10$ & $2^3 \cdot 3^2 \cdot 5^3$ & 4 & 2 & $x^4+32x^3+384x^2-512x+4096$\\
     \hline
 \end{tabular}
  \end{center}
\end{table}

In the previous table, the number fields defined by $f$ in the first and second row have the same Hilbert class field, and the same subfield satisfies the criteria of Test \ref{K1} in each case.

\subsection{Additional examples applying Test \ref{Km}}

Just as we used Test \ref{Km} to solve Problem \ref{CEP} in Examples 3.10 and 3.11, we apply Test \ref{Km} to the following examples. These are examples in which Test \ref{K1} failed, as noted in Table 2.

\begin{example}
Let $K=\bb{Q}(\alpha) \cong \bb{Q}[x]/(f)$ for $f=x^2+13$. The fractional ideal $I=2\bb{Z} \oplus (1+\alpha)\bb{Z}$ is non-principal with $N_{K/\bb{Q}}(I)=4$. We choose modulus $\mathfrak{m}=3\bb{Z}[\alpha]$ which is relatively prime to $I$. The ray class field of $K$ with modulus (and conductor) $\mathfrak{m}$ is the number field defined by 
\begin{align*}
&x^{16}  - 32x^{14} + 5676x^{12} - 316256x^{10} + 16999606x^8 - 689734368x^6 + \\
 &   19963313676x^4 - 378823819680x^2 + 3324557815569
 \end{align*}
over $\bb{Q}$. 

The subfield $F$ of $K_\mathfrak{m}$ defined by $x^2-74x+3721$ satisfies the criteria of Test \ref{Km}. Therefore, there is a single $E$-conjugacy class within $\mathcal{M}_f$ for $E=\mathcal{O}_F$.
\end{example}

\begin{example}
Let $f=x^2+14$, and $K=\bb{Q}(\alpha) \cong \bb{Q}[x]/(f)$. We consider the non-principal ideal $I=3\bb{Z} \oplus (-1+\alpha)\bb{Z}$. Since $N_{K/\bb{Q}}(I)=9$, we pick $\mathfrak{m}=2\bb{Z} \oplus \alpha \bb{Z}$, a prime ideal dividing $(2)$.

The ray class field $K_\mathfrak{m}$ is defined by 
\begin{align*}
&x^{16} - 32x^{14} + 5676x^{12} - 316256x^{10} + 16999606x^8 - 689734368x^6 + \\
    &19963313676x^4 - 378823819680x^2 + 3324557815569
    \end{align*}
     and has conductor $\mathfrak{m}$.
Using the criteria of Test \ref{Km}, we find the subfield $F$ of $K_\mathfrak{m}$ defined by 
\begin{align*}
x^4 - 1548x^3 - 3055050x^2 - 2822525676x + 3324557815569.
\end{align*}
Then for $E=\mathcal{O}_F$, all matrices in $\mathcal{M}_f$ are $E$-conjugate.
\end{example}

\vspace{3 in}

\newpage

\bibliographystyle{alpha}
\bibliography{conjugacy_extension}

\end{document}